\documentclass[11pt,a4paper]{article}
\usepackage[latin1,utf8]{inputenc}
\usepackage{amsmath}
\usepackage{amsfonts}
\usepackage{amssymb}

\usepackage[english]{babel}
\usepackage{amsmath,amsthm}
\usepackage{amsfonts}
\usepackage{appendix}
\usepackage{natbib}

\usepackage{graphicx,pict2e}
\usepackage{rotating}
\usepackage{algorithm}
\floatname{algorithm}{Algorithm}

\newtheorem{theo}{Theorem}
\newtheorem{prop}{Proposition}

\numberwithin{equation}{section}

\title{Large sample properties of the Midzuno sampling scheme with probabilities proportional to size}
\date{\today}

\author{Guillaume Chauvet\thanks{ENSAI/IRMAR, Campus de Ker Lann, 35170 Bruz, France. E-mail: chauvet@ensai.fr}}

\begin{document}

\maketitle

\begin{abstract}
\noindent Midzuno sampling enables to estimate ratios unbiasedly. We prove the asymptotic equivalence between Midzuno sampling and simple random sampling for the main statistical purposes of interest in a survey.
\end{abstract}


\normalsize

\section{Introduction} \label{sec:int}

\noindent \citet{mid:52} \citep[see also][]{sen:53} proposed a sampling algorithm to select a sample with unequal probabilities, while estimating unbiasedly a ratio. It may be of interest with a moderate sample size, when the small sample bias may be appreciable. Midzuno sampling has been recently considered in \citet{esc:ber:13} and \citet{hid:kim:nam:16}, for example. \\
\noindent We introduce a coupling algorithm between Midzuno sampling and simple random sampling, which enables to prove that the Horvitz-Thompson associated to these
two procedures are asymptotically equivalent. We obtain a central-limit theorem for the estimator of a total and for the estimator of a ratio. We also prove that
variance estimators suitable for simple random sampling are also consistent for Midzuno sampling. \\
\noindent The paper is organized as follows. The notation is introduced in Section \ref{sec:not}. The coupling procedure is described in Section
\ref{sec:cou}. It is used in Section \ref{sec:int} to prove the asymptotic normality of total and ratio estimators, and to establish the consistency of the proposed variance
estimators. Their behaviour is studied in Section \ref{sec:sim} through a simulation study with various sample sizes. We conclude in Section \ref{sec:conc}. The proofs are given in the Supplementary Material.

\section{Notation and assumptions} \label{sec:not}

\noindent We consider a finite population $U$ of size $N$, with a variable of interest $y$ taking the value $y_k$ for the unit $k \in U$. We are interested in estimating the total $Y = \sum_{k \in U} y_k$ or the ratio $R =Y/X$ with $X=\sum_{k \in U} x_k$ and $x_k>0$ is an auxiliary variable known for any unit $k \in U$. \\
\noindent Let $p_k>0$ be some probability for unit $k$, with $\sum_{k \in U} p_k = 1$. If the probabilities are chosen proportional to $x_k$, we have $p_k=x_k/X$. A sample $S$ of size $n$ is selected according to some sampling design with $\pi_k>0$ the inclusion probability of unit $k$. With Midzuno sampling, $p_k$ is the probability that the unit $k$ is selected at the first draw, while $\pi_k$ is the overall probability that the unit $k$ is selected in the sample, see Section \ref{ssec:not:mi}. The Horvitz-Thompson (HT) estimator for the total is $\hat{Y} = \sum_{k \in S} \frac{y_k}{\pi_k}$, and the substitution estimator for the ratio is $\hat{R}=\hat{Y}/\hat{X}$, with $\hat{X}=\sum_{k \in S} \frac{x_k}{\pi_k}$.

\subsection{Simple random sampling} \label{ssec:not:si}

\noindent If the sample is selected by simple random sampling in $U$, which is denoted as $SI(n;U)$, we obtain $\pi_k^{SI} = n/N$ and the estimators are
    \begin{eqnarray}  \label{sec:not:eq4}
      \hat{Y}_{SI} = \frac{N}{n} \sum_{k \in S_{SI}} y_k ~~~\textrm{and}~~~ \hat{R}_{SI} = \frac{\sum_{k \in S_{SI}} y_k}{\sum_{k \in S_{SI}} x_k}.
    \end{eqnarray}
The variance of the HT-estimator is
    \begin{eqnarray}  \label{sec:not:eq5}
      V(\hat{Y}_{SI}) = \frac{N(N-n)}{n} S_y^2 ~~~\textrm{with}~~~ S_y^2 = \frac{1}{N-1} \sum_{k \in U} \left(y_k-\frac{Y}{N}\right)^2,
    \end{eqnarray}
and is unbiasedly estimated by
    \begin{eqnarray}  \label{sec:not:eq6}
      \hat{V}(\hat{Y}_{SI}) = \frac{N(N-n)}{n} s_{y,SI}^2 ~~~\textrm{with}~~~ s_{y,SI}^2 = \frac{1}{n-1} \sum_{k \in S_{SI}} \left(y_k-\frac{\hat{Y}_{SI}}{N}\right)^2.
    \end{eqnarray}

\noindent Noting $z_k=y_k-R x_k$ and $\hat{z}_k=y_k-\hat{R} x_k$, the linearization variance approximation for $\hat{R}_{SI}$ is
    \begin{eqnarray}  \label{sec:not:eq7}
      V_{lin}(\hat{R}_{SI}) = \frac{N(N-n)}{n~X^2} S_z^2 ~~~\textrm{with}~~~ S_z^2 = \frac{1}{N-1} \sum_{k \in U} \left(z_k-\frac{\sum_{l \in U} z_l}{N}\right)^2,
    \end{eqnarray}
and the assorted variance estimator is
    \begin{eqnarray}  \label{sec:not:eq8}
      \hat{V}_{lin}(\hat{R}_{SI}) = \frac{N(N-n)}{n~\hat{X}_{SI}^2} s_{\hat{z},SI}^2 ~~~\textrm{with}~~~ s_{\hat{z},SI}^2 = \frac{1}{n-1} \sum_{k \in S_{SI}} \left(\hat{z}_k-\frac{\sum_{l \in S_{SI}} \hat{z}_l}{n}\right)^2.
    \end{eqnarray}
We prove in Section \ref{sec:int} that $\hat{V}$ and $\hat{V}_{lin}$ are consistent for Midzuno sampling.

\subsection{Midzuno sampling} \label{ssec:not:mi}

\noindent Suppose that the sample $S_{MI}$ is selected by means of the \citet{mid:52} sampling scheme, which is denoted as $MI$. A first unit ($k_1$, say) is selected in $U$ with probabilities $p_k$. A sample $S'_{MI}$ is then selected among the remaining units by $SI(n-1;U \setminus \{k_1\})$. The final Midzuno sample is $S_{MI}=S'_{MI} \cup \{k_1\}$, and the associated inclusion probabilities are
    \begin{eqnarray} \label{sec:not:eq7}
      \pi_k^{MI} & = & \frac{n-1}{N-1} + p_k \left(\frac{N-n}{N-1}\right).
    \end{eqnarray}
The main advantage of MI is that $\hat{R}_{MI}$ is exactly unbiased for $R$ if the probabilities $p_k$ are proportional to $x_k$.

\subsection{Assumptions} \label{ssec:not:ass}

\noindent We work under the asymptotic set-up of \citet{isa:ful:82}, where $U$ is embedded into a nested sequence of finite populations with $n,N \to \infty$. We suppose that the sampling rate is not degenerate, i.e. some constant $f \in ]0,1[$ exists s.t. $n/N \to f$. We will consider the following assumptions:
\begin{itemize}
  \item[H1:] Some constants $c_1,C_1$ exist, s.t. $0<c_1\leq N p_k \leq C_1$ for any $k \in U$.
  \item[H2:] Some constant $M$ exists, s.t. $N^{-1} \sum_{k \in U} y_k^4 \leq M$.
  \item[H3a:] Some constant $m_1>0$ exists, s.t. $S_y^2 \geq m_1$.
  \item[H3b:] Some constant $m_2>0$ exists, s.t. $S_z^2 \geq m_2$.
\end{itemize}

\section{Coupling procedure} \label{sec:cou}

\noindent The coupling procedure introduced in Algorithm \ref{algo:1} enables the justification of the closeness between MI and SI, as proved in Proposition \ref{prop2}.

\begin{algorithm}[t]
\begin{enumerate}
	\item Select some unit ($k_1$, say) in $U$ with probabilities $p_k$.
    \item Select $S'_{MI}$ by $SI(n-1;U \setminus \{k_1\})$. The MI sample is $S_{MI}=S'_{MI} \cup \{k_1\}$.
    \item Select some unit ($k_2$, say) in $U \setminus S'_{MI}$, with probability $n/N$ for $k_1$ and $1/N$ otherwise. The SI sample is $S_{SI}=S'_{MI} \cup \{k_2\}$.
 \end{enumerate}
\caption{Coupling procedure between MI and SI sampling} \label{algo:1}
\end{algorithm}

\begin{prop} \label{prop1}
  The sample $S_{SI}$ in Algorithm \ref{algo:1} is selected by $SI(n;U)$.
\end{prop}

\begin{prop} \label{prop2}
  Suppose that $S_{MI}$ and $S_{SI}$ are selected by Algorithm \ref{algo:1}, and that assumptions (H1)-(H2) hold. Then
    \begin{eqnarray} \label{prop2:eq1}
      E\left[\left(\hat{Y}_{MI}-\hat{Y}_{SI}\right)^4\right] = O(N^4 n^{-4})
      ~~ \textrm{and} ~~ E\left[\left(\hat{Y}_{MI}-Y\right)^4\right] = O(N^4 n^{-2}).
    \end{eqnarray}
\end{prop}

\noindent The first part of equation (\ref{prop2:eq1}) implies in particular that
    \begin{eqnarray} \label{sec:cou:eq2}
      \left(\sqrt{V(\hat{Y}_{MI})}-\sqrt{V(\hat{Y}_{SI})}\right)^2 & = & O(N^2 n^{-2}) = o\{V(\hat{Y}_{SI})\}.
    \end{eqnarray}
Consequently, $\hat{Y}_{MI}$ and $\hat{Y}_{SI}$ have asymptotically the same variance.

\section{Interval estimation} \label{sec:int}

\begin{theo} \label{theo1}
  Suppose that assumptions (H1), (H2) and (H3a) hold. Then
    \begin{eqnarray}
      \{V(\hat{Y}_{MI})\}^{-0.5} \{\hat{Y}_{MI}-Y\} & \longrightarrow_{\mathcal{L}} & \mathcal{N}(0,1), \label{theo1:eq1} \\
      E\left[ N^{-2} n \left\{\hat{V}(\hat{Y}_{MI})-V(\hat{Y}_{MI})\right\}\right]^2 & = & O(n^{-1}), \label{theo1:eq2}
    \end{eqnarray}
  with $\to_{\mathcal{L}}$ the convergence in distribution, and where $\hat{V}(\hat{Y}_{MI})$ is the SI variance estimator given in (\ref{sec:not:eq6}), applied to the sample $S_{MI}$.
\end{theo}

\noindent Theorem \ref{theo1} implies that the HT-estimator is asymptotically normally distributed under MI, and that the SI variance estimator is also consistent for MI, in the sense that $\{V(\hat{Y}_{MI})\}^{-1} \hat{V}(\hat{Y}_{MI}) \to_{Pr} 1$, with $\to_{Pr}$ the convergence in probability.
%
We now consider ratio estimation. We suppose that the $p_k$'s are defined proportionally to $x_k$, and we strengthen (H1) as
\begin{itemize}
  \item[H1b:] Some constants $c_1,C_1$ exist, s.t. $0<c_1\leq x_k \leq C_1$ for any $k \in U$.
\end{itemize}

\begin{prop} \label{prop3}
Suppose that assumptions (H1b) and (H2) hold. Then
    \begin{eqnarray}
      E\left[ \left\{(\hat{R}_{MI}-R)-X^{-1}(\hat{Z}_{MI}-Z) \right\}^2 \right] & = & O(n^{-2}).
    \end{eqnarray}
\end{prop}

\noindent This proposition entails in particular the validity of the linearization variance estimation, since $\{V_{lin}(\hat{R}_{SI})\}^{-1} V(\hat{R}_{SI}) \to 1$ if (H3b) is verified.

\begin{theo} \label{theo2}
  Suppose that assumptions (H1b), (H2) and (H3b) hold. Then
    \begin{eqnarray}
      \{V_{lin}(\hat{R}_{MI})\}^{-0.5} \{\hat{R}_{MI}-R\} & \longrightarrow_{\mathcal{L}} & \mathcal{N}(0,1), \label{theo2:eq1} \\
      E\left| n \left\{\hat{V}_{lin}(\hat{R}_{MI})-V_{lin}(\hat{R}_{MI})\right\}\right| & = & O(n^{-0.5}), \label{theo2:eq2}
    \end{eqnarray}
   where $\hat{V}_{lin}(\hat{R}_{MI})$ is the linearization SI variance estimator given in (\ref{sec:not:eq8}), applied to the sample $S_{MI}$.
\end{theo}

\noindent Theorem \ref{theo2} implies that the confidence interval $[\hat{R}_{MI} \pm u_{1-\alpha} \{\hat{V}_{lin}(\hat{R}_{MI})\}^{0.5}]$ has an asymptotic coverage
of $100(1-2\alpha) \% $.

\section{Simulation study} \label{sec:sim}

\noindent We conducted a simulation to evaluate the proposed variance estimators with small to moderate samples. We generated a population of $N=10,000$ units, with auxiliary variable $x$ generated according to a gamma distribution with shape and scale parameters $2$ and $5$, and we shifted and scaled the values so that $x_k$ lies between $1$ and $20$. We generated a variable of interest $y$ according to the model $y_k=x_k+\sigma~\epsilon_k$, with the $\epsilon_k$'s generated according to a standard normal distribution, and where
$\sigma$ was chosen so that the coefficient of determination was approximately $0.70$. \\
\noindent We repeated $B=20,000$ times MI, with $n$ ranging from $20$ to $500$. For each sample, the first unit $k_1$ is selected with probabilities $p_k$ proportional to $x_k$ by means of a fixed-size sampling algorithm, so that one unit exactly is selected. The $n-1$ other units of the MI sample are selected by simple random sampling in the rest of the population. To measure the bias of the estimator $\hat{\theta}$ of a parameter $\theta$, we used the Monte Carlo Percent Relative Bias
    \begin{eqnarray} \label{rel:bias}
      RB\{\hat{\theta}\} & = & 100 \times \frac{B^{-1} \sum_{b=1}^B \hat{\theta}_b-\theta}{\theta},
    \end{eqnarray}
where $\hat{\theta}_b$ denotes the estimator $\hat{\theta}$ in the b-th sample. We computed the relative bias for the estimators $\hat{Y}_{MI}$ and $\hat{R}_{MI}$. To measure the bias of some variance estimator $\hat{V}(\hat{\theta})$, we computed
    \begin{eqnarray} \label{rel:bias:2}
      RB\{\hat{V}(\hat{\theta})\} & = & 100 \times \frac{B^{-1} \sum_{b=1}^B \hat{V}_b(\hat{\theta}_b)-{MSE}(\hat{\theta})}{{MSE}(\hat{\theta})},
    \end{eqnarray}
where $\hat{V}_b(\hat{\theta}_b)$ denotes the variance estimator in the b-th sample, and where ${MSE}(\hat{\theta})$ is a simulation-based approximation of the true mean square error obtained from an independent run of $100,000$ simulations. As a measure of stability of $\hat{V}(\hat{\theta})$, we used the Relative Root Mean Square Error
    \begin{eqnarray*} \label{instab}
      RRMSE\{\hat{V}(\hat{\theta})\} & = & 100 \times \frac{\left[B^{-1}\sum_{b=1}^B \left\{\hat{V}_b(\hat{\theta}_b)-{MSE}(\hat{\theta}) \right\}^2\right]^{1/2}}{{MSE}(\hat{\theta})}.
    \end{eqnarray*}
We computed the relative bias and the relative root mean square error for the variance estimators $\hat{V}(\hat{Y}_{MI})$ and $\hat{V}_{lin}(\hat{R}_{MI})$. Finally, we computed the error rate of the normality-based confidence intervals with nominal one-tailed error rate of 2.5 \% in each tail. \\

\noindent The results are given in Table \ref{tab:sim}. We first consider $\hat{Y}_{MI}$, which is always unbiased as expected. The estimator $\hat{V}(\hat{Y}_{MI})$ is positively biased with small sample sizes, but the bias vanishes when $n$ grows. Despite the variance being overestimated, the coverage rates are well respected in any case and are even below the nominal level for small sample sizes. This is likely due to the fact that the asymptotic normality is a crude approximation when $n$ is small, and that the Student $t$-distribution would presumably perform better. With $n=20$, the $2.5 \% $ quantile of the $t$-distribution with $n-1=19$ degrees of freedom is $u_{0.025}^{Stu}=2.093$. Using the $2.5 \% $ normal quantile $u_{0.025}^{Nor}=1.96$ instead therefore leads to narrowing the confidence interval, which compensates for overestimating the variance. As for the RRMSE, we note that it decreases when $n$ grows, as expected. We now turn to $\hat{R}_{MI}$. It is unbiased in all the cases considered, as expected. The estimator $\hat{V}(\hat{Y}_{MI})$ is almost unbiased and the coverage rates are well respected in all cases.

\begin{sidewaystable}[h!]
\caption{Relative bias of point estimators, Relative Bias and Relative Root Mean Square Error of variance estimators, and coverage rates} \label{tab:sim}
	\begin{tabular}{|l|c|ccc|c|ccc|} \hline
            & $\hat{Y}_{MI}$ & \multicolumn{3}{|c|}{$\hat{V}(\hat{Y}_{MI})$} & $\hat{R}_{MI}$ & \multicolumn{3}{|c|}{$\hat{V}_{lin}(\hat{R}_{MI})$} \\ \hline
            & RB (\%) & RB (\%) & RRMSE & Cov. Rate & RB (\%) & RB (\%) & RRMSE & Cov. Rate \\ \hline
    $n=20$  & 0.0 & 12.8 & 51.8 & 94.3 & 0.0 & 0.7 & 40.5 & 94.1 \\
    $n=40$  & 0.0 & 7.2  & 34.0 & 94.8 & 0.0 & 0.9 & 28.6 & 94.6 \\
    $n=60$  & 0.0 & 5.0  & 26.8 & 94.9 & 0.0 & 0.7 & 23.3 & 94.7 \\
    $n=80$  & 0.0 & 3.9  & 22.9 & 94.7 & 0.0 & 0.6 & 20.0 & 94.8 \\
    $n=100$ & 0.0 & 2.8  & 20.2 & 94.8 & 0.0 & 0.4 & 18.0 & 94.9 \\
    $n=200$ & 0.0 & 1.7  & 14.1 & 94.9 & 0.0 & 0.3 & 12.6 & 94.9 \\
    $n=500$ & 0.0 & 0.1  & 8.5  & 95.1 & 0.0 & 0.2 & 7.8  & 95.2 \\ \hline
	\end{tabular}
\end{sidewaystable}

\section{Conclusion} \label{sec:conc}

\noindent In this paper, we have proved rigorously that Midzuno sampling is equivalent to simple random sampling for main statistical purposes. This is also justified empirically by the simulation results, with a small sample size $n=20$. Despite the large number of papers which have considered this method ($275$ according to GoogleScholar), it seems therefore of limited interest. \\
\noindent From equation (\ref{sec:not:eq7}), the range of possible inclusion probabilities under Midzuno sampling is very limited. \cite{dev:til:98} have proposed a generalization of the Midzuno method, suitable for any set of inclusion probabilities. Extending the results of the current paper to the generalized Midzuno method would be an interesting matter for further research.

\bibliographystyle{apalike}

\end{document}